\documentclass{amsart}
\usepackage{amscd, amsmath, amsthm, amssymb}

\newtheorem{theorem}{Theorem}[section]
\newtheorem{lemma}[theorem]{Lemma}
\newtheorem{proposition}[theorem]{Proposition}

\newtheorem{corollary}[theorem]{Corollary}

\theoremstyle{definition}     
\newtheorem{definition}[theorem]{Definition}

\newtheorem{claim}[theorem]{Claim}
\newtheorem{question}[theorem]{Question}

\theoremstyle{remark}
\newtheorem{remark}[theorem]{Remark}

\numberwithin{equation}{section}

\begin{document}

\title[hyperk\"aher manifolds ]
{Automorphisms of hyperk\"ahler manifolds in the view of topological entropy}

\author[K. Oguiso]{Keiji Oguiso}
\address{Graduate School of Mathematical Sciences,
University of Tokyo, Komaba, Meguro-ku,
Tokyo 153-8914, Japan, and Korea Institute for Advanced Study, 
207-43 Cheonryangni-2dong, Dongdaemun-gu, Seoul 130-722, Korea}
\email{oguiso@ms.u-tokyo.ac.jp}

\dedicatory{Dedicated to Professor Igor Dolgachev on the occasion of
his sixtieth birthday}


\begin{abstract} 
First we show that any group of automorphisms of null-entropy 
of a projective hyperk\"ahler manifold $M$ is almost abelian of rank at most 
$\rho(M) - 2$.  
We then characterize automorphisms of a K3 surface with null-entropy and those with positive entropy 
in algebro-geometric terms. 
We also give an example of a group of automorphisms which is not almost abelian in each dimension. 
\end{abstract}

\maketitle


\setcounter{section}{0}
\section{Introduction - Background and main results} 
The aim of this note is to study groups of automorphisms of 
a hyperk\"ahler manifold from two points of view: topological entropy and how close to (or how far from) 
abelian groups. Our main results are Theorems (1.3) and (2.1). 
\par \vskip 1pc  
\noindent 
{\bf 1.} Let $M$ be a compact K\"ahler manifold. We denote the biholomorphic 
automorphism group of $M$ by ${\rm Aut}\,(M)$. By the fundamental work of Yomdin, 
Gromov and Friedland, 
the {\it topological entropy} $e(g)$ of an automorphism 
$g \in \text{Aut}\,(M)$ can be computed in three different ways: topological, differential-geometrical, 
and cohomological (see [Yo], [Gr], [Fr], also [DS2]). In this note, we employ the cohomological one as its definition:   
$$e(g)\,\, := \log \delta(g)\,\, .$$ 
Here $\delta(g)$ is the spectral radius of the action of $g$ on the cohomology ring $H^{*}(M, \mathbf C)$, i.e. 
the maximum of the absolute values of eigenvalues of 
the $\mathbf C$-linear extension of $g^{*} \vert H^{*}(M, \mathbf Z)$. 
One has $e(g) \geq 0$, and $e(g) = 0$ iff the eigenvalues of 
$g^{*}$ are on the unit circle $S^{1} := \{\, z \in \mathbf C\, \vert\, \vert z \vert = 1\,\}$. 
Furthermore, by Dinh and Sibony 
[DS1], $e(g) > 0$ iff 
some eigenvalues of $g^{*} \vert H^{1,1}(M)$ are {\it outside} the unit circle 
$S^{1}$. A subgroup $G$ of 
$\text{Aut}\, (M)$ is said to be of {\it null-entropy} (resp. of 
{\it positive entropy}) if $e(g) = 0$ 
for $\forall g \in G$ (resp. $e(g) > 0$ for $\exists g \in G$). 
\par \vskip 1pc  
\noindent 
{\bf 2.} Next we recall a few facts about hyperk\"ahler manifolds. K3 surfaces are nothing but $2$-dimensinal 
hyperk\"ahler manifolds. All what we need is 
reviewed in [Og2, Section 2]:  

\begin{definition} \label{definition:hk} {\it A hyperk\"ahler manifold} 
is a compact complex simply-connected K\"ahler 
manifold $M$ admitting an everywhere non-degenerate global 
holomorphic $2$-form $\omega_{M}$ such that $H^{0}(M, \Omega_{M}^{2}) = 
\mathbf C \omega_{M}$. 
\end{definition} 

Let $M$ be a hyperk\"ahler manifold. Then the second cohomology group 
$H^{2}(M, \mathbf Z)$ admits a natural $\mathbf Z$-valued 
symmetric bilinear form of signature $(3, 0, b_{2}(M) -3)$, called Beauville-Bogomolov-Fujiki's form 
(BF-form for short). 
BF-form is exactly the cup-product when $M$ is a K3 surface. 
The signature of the N\'eron-Severi group $NS(M)$ w.r.t. BF-form is either 
$(1, 0, \rho(M) -1)$, $(0, 1, \rho(M) -1)$, or $(0, 0, \rho(M))$. 
Here $\rho(M)$ is the Picard number of $M$. We call these three cases 
{\it hyperbolic}, {\it parabolic} and {\it elliptic} respectively. 
Due to Huybrechts [Hu], $M$ is 
projective iff $NS(M)$ is hyperbolic. Note also that, in dimension $2$, 
$NS(M)$ is hyperbolic, parabolic, elliptic iff the algebraic dimension 
$a(M)$ is $2$, $1$, $0$, respectively (see eg. [BPV]). 
\par \vskip 1pc  

In our previous note, we have shown the following:

\begin{theorem} \label{theorem:prev} {\rm [Og2]} The bimeromorphic automorphism group ${\rm Bir}\, (M)$ 
of a non-projective hyperk\"ahler manifold $M$ is an almost abelian group of finite rank. 
More precisely, ${\rm Bir}\, (M)$ is an almost abelian group of rank at most $1$ (resp. $\rho(M) -1$) 
when $NS(M)$ is elliptic (resp. parabolic). 
\end{theorem} 

We call a group $G$ {\it almost 
abelian} (resp. {\it almost abelian of finite rank} $r$) 
if there are a normal subgroup $G^{(0)}$ of $G$ of finite index, 
a finite group $K$ and an abelian group $A$ (resp. $A = \mathbf Z^{r}$)  which fit in 
the exact sequence
$$1 \longrightarrow K \longrightarrow G^{(0)} \longrightarrow A \longrightarrow 0\,\,  .$$

The rank $r$ is well-defined, and invariant under 
replacing $G$ by a subgroup $H$ of finite index 
and by a quotient group $Q$ of $G$ by a finite normal subgroup (cf. [Og2, Section 9]). 
It is then clear that a subgroup of an almost abelian group of rank $r$ is an almost abelian group 
of rank at most $r$.  
\par \vskip 1pc  
\noindent 
{\bf 3.} Our main result is the following: 

\begin{theorem} \label{theorem:main} Let $M$ be a projective hyperk\"ahler manifold. Let $G < {\rm Aut}\, (M)$ 
(resp. $G < {\rm Bir}\,(M)$). Assume that $G$ is of null-entropy (resp. of null-entropy 
at $H^{2}$-level, that is, the eigenvalues of $g^{*} \vert H^{2}(M, \mathbf Z)$ are on the unit circle $S^{1}$). 
Then $G$ is an almost abelian group of rank at most $\rho(M) -2$. Moreover, 
this estimate is optimal in ${\rm dim}\, M = 2$. 
\end{theorem}

The key ingredient is Theorem (2.1), a result of 
linear 
algebra, whose source has been back to an important observation of 
Burnside [Bu1], and Lie's Theorem (cf. [Hm]). We shall prove Theorem (1.3) in Section 3. 
\par \vskip 1pc  
\noindent 
{\bf 4.} Next, as an application of Theorem (1.3), we shall 
reproduce the following algebro-geometric characterization of positivity 
of entropy of automorphisms of a K3 surface. We should notice that this result 
is essentially known and can be read from works of Cantat [Ca1, 2]:
 
\begin{theorem} \label{theorem:zar} 
Let $M$ be a (not necessarily projective) K3 surface, $G < {\rm Aut}\, (M)$, and 
$g \in {\rm Aut}\, (M)$. Then: 
\begin{list}{}{
\setlength{\leftmargin}{10pt}
\setlength{\labelwidth}{6pt}
}
\item[(1)] $G$ is of null-entropy iff either $G$ is finite or $G$ makes an 
elliptic fibration $\varphi : M \longrightarrow \mathbf P^{1}$ stable. 
In particular, if $G$ is of null-entropy, then $G$ has no Zariski dense orbit, 
and conversely, if $G$ makes an elliptic fibration stable, then $G$ is almost abelian 
of finite rank. 
\item[(2)] $e(g) > 0$ iff $g$ has a Zariski dense orbit, i.e. there 
is a point $x \in M$ such that the set $\{g^{n}(x) \vert n \in \mathbf Z\}$ 
is Zariski dense in $M$. 
\end{list}
\end{theorem} 

Theorem (1.4) is proved in Sections 3 and 4. This theorem is also motivated by earlier observations of 
[Sn] and by the following question posed 
by McMullen: 

\begin{question} \label{question:Mc} {\rm [Mc]} Does a K3 automorphism $g$ 
have a dense orbit (in the Euclidean topology) 
when a K3 surface is projective and $e(g) > 0$?
\end{question} 

In the same paper [Mc], he constructed a K3 surface $M$ of $\rho(M) = 0$ having an automorphism $g$ with Siegel disk. 
For this $(M, g)$, one knows that $e(g) > 0$, $g$ has a Zariski dense orbit but no dense orbit in the Euclidean 
topology, and that ${\rm Aut}\,(M) \simeq \mathbf Z$ (See [ibid] and also [Og2]).
\par \vskip 1pc  
\noindent 
{\bf 5.} So far, all groups in consideration are almost abelian. As a sort of counter parts, we shall show 
the following Theorem in Section 5:
\begin{theorem} \label{theorem:nonab} 
\begin{list}{}{
\setlength{\leftmargin}{10pt}
\setlength{\labelwidth}{6pt}
}

\item[(1)] Let $M$ be a K3 surface admitting two different Jacobian fibrations 
of positive Mordell-Weil rank. Then ${\rm Aut}\,(M)$ is not almost abelian (and hence is of positive entropy). 
This happens, for instance, for a K3 surface $M$ of maximum Picard number $\rho(M) = 20$. 
\item[(2)] In each dimension $2m$, there is a projective hyperk\"ahler manifold $M$ 
whose ${\rm Aut}\,(M)$ is not almost abelian (and hence is of positive entropy). 
\end{list}
\end{theorem} 

Theorem (1.6)(1) is a part of refinement of a result [Ca2] (see also [CF] for 
a generalization in foliated case) 
and also a slight generalization of a result of Shioda and Inose [SI]: 
{\it ${\rm Aut}\,(M)$ 
is an infinite group for a K3 surface $M$ with $\rho(M) = 20$.} In [ibid], this has been shown by finding a Jacobian fibration of positive Mordell-Weil rank. 
\par \vskip 1pc  
\noindent 

\par \vskip 1pc  
{\it Acknowledgement.} I would like to express my thanks to Professors 
S. Cantat, D. Huybrechts, Y. Ishii, S. Kawaguchi, Y. Kawamata, J.H. Keum, T. Shioda 
and De-Qi Zhang for their warm encouragements and valuable suggestions. 
The semi-final version has been completed during my stay at KIAS March 2005. 
I would like to express my thanks to Professors J.M. Hwang and B. Kim for invitation. 

\section{Group of isometries of null-entropy of a hyperbolic lattice}

The main result of this section is Theorem (2.1).

A lattice is a pair $L := (L, b)$ 
of a free abelian group $L \simeq \mathbf Z^{r}$ 
and a non-degenerate symmetric bilinear form 
$$b : L \times L \longrightarrow \mathbf Z\,\, .$$

A submodule $M$ (resp. an element $v \not= 0$) of $L$ is 
{\it primitive} if and only if $L/M$ 
(resp. $L/\mathbf Z v$)  is free.

A scalar extension $(L \otimes K, b \otimes id_{K})$ of $(L, b)$ by $K$ is 
written as $(L_{K}, b_{K})$. We often write $b_{K}(x, y)$, $b_{K}(x, x)$ ($x, y \in L_{K}$) 
simply as $(x, y)$, $(x^{2})$.

The non-negative integer $r$ is called the rank of $L$. The signature of $L$ 
is the signature, i.e. the numbers of positive-, zero-, negative-eigenvalues,  
of a symmetric matrix associated to $b_{\mathbf R}$. It is denoted 
by $\text{sgn}\, L$. 
The lattice $L$ is called {\it hyperbolic} (resp. {\it parabolic}, {\it elliptic}) 
if $\text{sgn}\, L$ is $(1, 0, r-1)$ (resp. $(0, 1, r-1)$, $(0, 0, r)$).

{\it In what follows, $L$ is assumed to be a hyperbolic lattice.}

The positive cone $\mathcal P(L)$ of $L$ is one of the two 
connected components of:

$$\mathcal P'(L_{\mathbf R}) := \{x \in L_{\mathbf R}\,\, 
\vert\,\, (x, x) > 0\}\,\, .$$ 

The boundary (resp. the closure) of $\mathcal P(L)$ is denoted by 
$\partial \mathcal P(L)$ (resp. $\overline{\mathcal P}(L)$). Obviously, $(x^{2}) = 0$ 
if $x \in \partial \mathcal P(L)$. 
Let $x, x' \in \overline{\mathcal P}(L) \setminus \{0\}$. Then, by the Schwartz inequality, 
$(x, x') \geq 0$ and the equality holds exactly when $x$ and $x'$ are proportional boundary points.

We denote the group of isometries of $L$ by:

$${\rm O}(L)\,\, := \,\, 
\{g \in {\rm Isom}_{\rm {group}}(L)\,\, \vert\,\, (gx, gy) = (x, y)\,\, 
\forall x, y \in L \}\,\, .$$ 

We have an index $2$ subgroup:

$${\rm O}(L)' := \{g \in {\rm O}(L)\,\, \vert\,\, 
g(\mathcal P(L)) = \mathcal P(L) \}\,\, .$$

Let $g \in {\rm O}(L)'$. The spectral value $\delta(g)$ of $g$ is 
the maximum of the absolute values of eigenvalues of $g_{\mathbf C}$. 
By abuse of notation, 
we call 
$$e(g) := \text{log}\, \delta(g)\,\,$$ 
the {\it entropy} of $g$. As we shall see in Proposition (2.2), 
$e(g) \geq 0$, and $e(g) = 0$ if and only if  
the eigenvalues of $g_{\mathbf C}$ lie on the unit circle 
$S^{1} := \{z \in \mathbf C\,\, \vert\,\, \vert z \vert = 1\}$. An element 
$g$ of ${\rm GL}\,(r, \mathbf C)$ is called {\it unipotent} if all the eigenvalues are $1$.  

The aim of this section is to prove the following: 

\begin{theorem} \label{theorem:mainlttice} Let $L$ be a hyperbolic lattice of rank $r$ and 
$G$ be a subgroup of ${\rm O}(L)'$. 
Assume that $G$ is an infinite group of null-entropy, i.e. 
$\vert G \vert = \infty$ and that $e(g) = 0$ for all $g \in G$.  
Then: 

\begin{list}{}{
\setlength{\leftmargin}{10pt}
\setlength{\labelwidth}{6pt}
}
\item[(1)] There is $v \in \partial \mathcal P(L) \setminus \{0\}$ such that 
$g(v) = v$ for all $g \in G$. 
Moreover, the ray $\mathbf R_{> 0} v$ in $\partial \mathcal P(L)$ is 
unique and defined over 
$\mathbf Z$. In other words, one can take unique such $v$ which is primitive 
in $L$. 
\item[(2)] There is a normal subgroup $G^{(0)}$ of $G$ such that 
$[G :G^{(0)}] < \infty$ and $G^{(0)}$ is isomorphic to a free abelian group of 
rank at most 
$r-2$. In particular, $G$ is almost abelian of rank at most $r-2$.  
\end{list}
\end{theorem} 
 
We shall prove Theorem (2.1) dividing into several steps. In what follows $G$ is as in Theorem (2.1).

\begin{proposition} \label{proposition:infty} 

\begin{list}{}{
\setlength{\leftmargin}{10pt}
\setlength{\labelwidth}{6pt}
}
\item[(1)] The eigenvalues of $g \in G$ lie on the unit circle $S^{1}$. 
\item[(2)] There is a positive integer $n$ such that $g^{n}$ is 
unipotent for all $g \in G$. 
\item[(3)] There is $g \in G$ such that ${\rm ord}\, g = \infty$. 
\end{list}
\end{proposition} 

\begin{proof} 
Since $g \in {\rm O}(L)$ and $L$ is non-degenerate, we have 
$\text{det}(g) = \pm 1$. 
Let $\alpha_{i}$ ($1 \leq i \leq r$) be the eigenvalues of 
$g_{\mathbf C}$ (counted with multiplicities). 
Then 
$$1\,\, = \vert {\rm det}\, g \vert \,\, = \,\, \Pi_{i=1}^{r} \vert \alpha_{i} \vert\,\, .$$ 
Thus,  
$e(g) \geq 0$, and $e(g) = 0$ if and only if  
the eigenvalues of $g_{\mathbf C}$ lie on $S^{1}$. 
This proves (1). 

Since $g$ is 
defined over $\mathbf Z$, the eigenvalues of $g_{\mathbf C}$ are 
all algebraic integers. 
Thus they are all roots of unity by (1) and by the Kronecker Theorem [Ta]. 
Since the eigen polynomial of $g$ is of degree $r$ and it is now the product 
of cyclotomic polynomials, the eigenvalues of $g$ lie on:  
$$\cup_{\varphi(d) \leq r}\,\, \mu_{d}\,\, < S^{1}\,\, .$$ 
Here $\mu_{d} = \langle \zeta_{d} \rangle 
:= \{z \in \mathbf C\,\, \vert\,\, z^{d} = 1\}$ and 
$\varphi(d) := \vert \text{Gal}(\mathbf Q(\zeta_{d})/\mathbf Q) \vert$ is 
the Euler function.
There are finitely many $d$ with $\varphi(d) \leq r$. 
Let $n$ be their product. Then $g^{n}$ is unipotent for all $g \in G$. 
This proves (2).

Let us show (3).  We have $G < \text{GL}(L_{\mathbf C}) \simeq \text{GL}(r, \mathbf C)$. 
Suppose to the contrary that $G$ consists of elements of finite order. 
Then each $g \in G$ is diagonalizable over $L_{\mathbf C}$. Let 
$n$ be as in (2). 
Then, $g^{n} = id$ for all $g \in G$, i.e. $G$ has a finite exponent $n$. 
However, then $\vert G \vert < \infty$ by the following theorem due to 
Burnside [Bu2] (based on [Bu1]), a contradiction: 

\begin{theorem} \label{theorem:burnside} Any subgroup of 
${\rm GL}\,(r, \mathbf C)$ 
of finite exponent is finite, i.e. $\vert G \vert < \infty$. 
\end{theorem} 
\end{proof}

\begin{remark} \label{remark:burnsidepb} 
After Theorem (2.3), Burnside asked if a group of finite exponent is finite 
or not (the Burnside problem). It is now known to be false in general 
even for finitely generated groups. The original Burnside problem 
has been now properly modified and completely solved by Zelmanov. 
(See for instance [Za] about Burnside problems.)
\end{remark} 

Set 
$$H := \{g \in G\,\, \vert\,\, g\,\, \text{is unipotent}\}\,\, .$$ 

\begin{lemma} \label{lemma:uni} 
$H$ is a normal subgroup of $G$. Moreover, $H$ has an element 
of infinite order. 
\end{lemma}

\begin{proof} Let $g$ be as in Proposition (2.2)(3) and $n$ be as in 
Proposition (2.2)(1). Then 
$g^{n}$ is an element of $H$ of infinite order. The only non-trivial 
part is the closedness of $H$ under the product, 
that is, if $a, b \in H$ then $ab \in H$. Since 
$a^{l}b^{m} \in G$, this follows from 
the next slightly more general result together with Proposition (2.2)(1). 
\end{proof} 

\begin{proposition} \label{proposition:quasiuni} 
\begin{list}{}{
\setlength{\leftmargin}{10pt}
\setlength{\labelwidth}{6pt}
}
\item[(1)] Let $z_{i} \in S^{1}$ ($1 \leq i \leq r$). Then 
$$\vert \sum_{i=1}^{r} z_{i} \vert \leq r\,\, ,$$ 
and 
$$\sum_{i=1}^{r} z_{i} = r\,\, \Longleftrightarrow\,\, z_{i} = 1\,\, 
\forall i\,\, .$$ 
\item[(2)] Let $A, B\, \in\, {\rm GL}\,(r, \mathbf C)$ such that 
$A$ and $B$ are 
unipotent 
and such that the eigenvalues of $A^{l}B^{m}$ lie on $S^{1}$ for all 
$l, m \in \mathbf Z_{\geq 0}$. 
Then $A^{l}B^{m}$ is unipotent for all  $l, m \in \mathbf Z_{\geq 0}$.
\end{list}
\end{proposition} 

\begin{proof}  The statement (1) follows from the triangle inequality. 
Let us show (2). We may assume that $A$ is of the Jordan form 
$$A := J(r_{1}, 1) \oplus \cdots \oplus J(r_{k}, 1)\,\, .$$ 
We fix $m$ and decompose $B^{m}$ into blocks as:  
$$B^{m} = (B_{ij})\,\, , \,\, B_{ij} \in \text{M}(r_{i}, r_{j}, \mathbf C)\,\, .$$

We have 

$$\text{tr}\, A^{l}B^{m} = \sum_{i=1}^{k} \text{tr}J(r_{i}, 1)^{l}B_{ii}\,\, .$$ 

Using an explicit form of $J(r_{i}, 1)^{l}$ and 
calculating its product with $B_{ii}$, we obtain:  

$$\text{tr}\, A^{l}B^{m} = b_{s}l^{s} + b_{s-1}l^{s-1} + \cdots + b_{1}l + 
\text{tr}\, B^{m}\,\, .$$ 
Here $s = \text{max}\, \{r_{i}\} - 1$ and $b_{t}$ are constant 
being independent of $l$. 

On the other hand, by the assumption and (1), we have 

$$\vert \text{tr}\, A^{l}B^{m} \vert \leq r\,\, ,$$ 

for all $l$. Thus, by varying $l$ larger and larger, we have 
$$b_{t} = 0\,\, \text{for}\,\, t = s\, ,\, s-1\, ,\, \cdots\, , 1$$ 
inductively. This implies 

$$\text{tr}\, A^{l}B^{m} = \text{tr}\, B^{m}\,\, ,$$ 
for all $l$. Since $B$ is unipotent, so is $B^{m}$. Thus,  
$\text{tr}\, B^{m} = r$ and therefore  
$$\text{tr}\, A^{l}B^{m} = r\,\, .$$
Since the eigenvalues of $A^{l}B^{m}$ lie on $S^{1}$, this implies the result. 
\end{proof} 

\begin{corollary} \label{corollary:upper} 
$H < T(r, \mathbf Q)$ under suitable basis $\langle u_{i} \rangle_{i=1}^{r}$ 
of $L_{\mathbf Q}$. 
Here $T(r, \mathbf Q)$ is the subgroup of ${\rm GL}\,(r, \mathbf Q)$ 
consisting of the upper triangle unipotent matrices.  
\end{corollary} 

\begin{proof} By Lemma (2.5), $H$ is a unipotent subgroup of 
$\text{GL}(r, \mathbf Q)$. Thus, the result follows from Lie's Theorem 
(see for instance [Hm, 17.5 Theorem]). Here is a small remark. Lie's Theorem 
in [Hm] is formulated and proved over algebraically closed field of characteristic $0$. 
So, precisely speaking, we have that $H < T(r, \mathbf C)$ under suitable basis 
$\langle u_{i} \rangle_{i=1}^{r}$ of 
$L_{\mathbf C}$. On the other hand, $L_{\mathbf Q}$ and all the elements of $H$ are defined over $\mathbf Q$. 
Thus, such basis $\langle u_{i} \rangle_{i=1}^{r}$ of $L_{\mathbf C}$ are nothing but solutions of a system of 
linear equations with rational coefficients in the range ${\rm det}\, (u_{i}) \not= 0$. Thus, exisitence of such basis 
over $\mathbf C$ implies the exisitence of the desired basis over $\mathbf Q$ as claimed.  
\end{proof}

Let $G$ and $H$ be as above. So far, we did not yet use the fact that 
$G$ and $H$ are subgroups of 
${\rm O}(L)$. From now, we shall use this fact.

\begin{lemma} \label{lemma:nulvect} 
\begin{list}{}{
\setlength{\leftmargin}{10pt}
\setlength{\labelwidth}{6pt}
}
\item[(1)] Let $N$ be a subgroup of ${\rm O}(L)$. 
Assume that there is an element $x$ of $L$ such that $(x^{2}) > 0$ and 
$a(x) = x$ for all $a \in N$. Then $\vert N \vert < \infty$. In particular, 
there is no such $x$ for $H$. 
\item[(2)] There is a unique ray $\mathbf R_{> 0} v (\not= 0)$ in 
$\partial \mathcal P(L)$ such that 
$a(v) = v$ for all $a \in H$. 
Moreover, the ray $\mathbf R_{> 0} v$ is defined over 
$\mathbf Z$, i.e. one can take unique such $v$ which is primitive in $L$. 
\item[(3)] Let $v$ be as in (2). Then $b(v) = v$ for all $b \in G$. 
\end{list}
\end{lemma} 
 
\begin{proof}  Let us first show (1). By assumption, $N$ can be naturally 
embedded into ${\rm O}(x^{\perp}_{L})$. Here $x^{\perp}_{L}$ is the orthogonal 
complement of $x$ in $L$. Since $L$ is hyperbolic and $(x^{2}) > 0$, 
the lattice $x^{\perp}_{L}$ is of negative definite.
Thus ${\rm O}(x^{\perp}_{L})$ is finite. This shows (1).

Let us show (2). By (1) and by $\vert H \vert = \infty$, we have $r \ge 2$. First, we find 
$v \in \partial \mathcal P(L) \cap L 
\setminus \{0\}$ 
such that $a(v) = v$ for all $a \in H$ by the induction on $r \geq 2$.

By Corollary (2.7), there is $u \in L \setminus \{0\}$ such that $a(u) = u$ for all $a \in H$. 
$H$ is then naturally embedded into ${\rm O}(u^{\perp}_{L})$ 
if $(u^{2}) \not= 0$. 

By (1), we have $(u^{2}) \leq 0$. If $(u^{2}) = 0$, then we are done. 
If $(u^{2}) < 0$, then $u^{\perp}_{L}$ is of signature $(1, 0, r-2)$. 
If in addition $r = 2$, then $u^{\perp}_{L}$ is of positive definite and 
${\rm O}(u^{\perp}_{L})$ is finite, a contradiction. 
Hence $(v^{2}) = 0$ when $r = 2$ and we are done when $r = 2$. If $r > 2$, 
then by the induction, we can find a desired $v$ in $u^{\perp}_{L}$ 
and we are done. 

If necessarily, by replacing $v$ by $\pm v/m$, we have that
$\mathbf R_{> 0} v \subset \partial \mathcal P(L)$ and $v \in L$ is primitive 
as well. 

Next, we shall show the uniqueness of $\mathbf R_{>0}v$. Suppose to the contrary that there is 
a ray $(0 \not= ) \mathbf R_{>0}u \subset \partial \mathcal P(L)$ s.t.  
$a(u) = u$ for 
all $a \in H$ and $\mathbf R_{>0} u \not= \mathbf R_{> 0} v$. Then 
the common eigenspace of eigenvalue $1$
$$V\, :=\, \cap_{h \in H} V(h, 1)$$ 
is a linear subspace of $L_{\mathbf C}$ of ${\rm dim}\, V \geq 2$. Since each 
$V(h, 1)$ is defined over $\mathbf Q$, so is $V$. Thus, we may find such $u$ 
in $L$. So, from the first, we may assume that $u \in L$. 
We have $((v + u)^{2}) > 0$ by the Schwartz inequality. 
We have also $a(v + u) = v + u$ for all $a \in H$. 
However, we would then have 
$\vert H \vert < \infty$ by (1), a contradiction. 
Now the proof of (2) 
is completed.

Finally, we shall show (3). Let $b \in G$. Put $b(v) = u$. Let $a \in H$. 
Since $H$ is a normal subgroup of $G$, 
there is $a' \in H$ such that $ab = ba'$. Then 
$$a(u) = a(b(v)) = b(a'(v)) = b(v) = u\,\, .$$ 
Verying $a$ in $H$ and using (2), we find that $u \in \mathbf R_{> 0}v$, i.e. $b(v) = \alpha v$ for 
some $\alpha > 0$. Since the eigenvalues of $b$ are on $S^{1}$, 
we have $\alpha = 1$. 

\end{proof}

\begin{proposition} \label{proposition:almstabel} 
Let $L$ be a hyperbolic lattice of rank $r$ and $N$ be a subgroup of ${\rm O}(L)$. 
Assume that there is a primitive element 
$v \in \partial \mathcal P(L) \cap L \setminus \{0\}$ 
such that 
$h(v) = v$ for all $h \in N$. Let  
$$\overline{L} := v^{\perp}_{L}/\mathbf Z v\,\, .$$ 
Then $\overline{L}$ is elliptic, of rank $r-2$, and 
the isometry $N$ on $L$ naturally descends to the isometry of $\overline{L}$, 
say $h \mapsto \overline{h}$. 
Set 
$$N^{0} := {\rm Ker}\, (N \longrightarrow {\rm O}(\overline{L}) 
\times \{\pm 1\}\,\, ;\,\, h\, \mapsto\, 
(\overline{h}\, ,\, {\rm det}\, h)\,)\,\, .$$ 
Then $N^{0}$ is of finite index in $N$ and $N^{0}$ 
is a free abelian group of rank at most $r-2$. Moreover $N$ is 
of null-entropy.

\end{proposition}

\begin{proof} The first part of the proposition is clear. 
We shall show the last two assertions. Since $\overline{L}$ is elliptic, the group 
${\rm O}(\overline{L})$ is finite. Hence 
$[N : N^{0}] \leq 2\cdot \vert {\rm O}(\overline{L}) \vert < \infty$. 
Choose an integral basis $\langle \overline{u}_{i} \rangle_{i=1}^{r-2}$ of $\overline{L}$. 
Let $u_{i} \in v^{\perp}_{L}$ be a lift of $\overline{u}_{i}$. Then 
$$\langle v, u_{i} \rangle_{i=1}^{r-1}$$ 
forms an integral basis of $v^{\perp}_{L}$. Since 
$v^{\perp}_{L}$ is primitive in $L$, there is an element $w \in L$ such that 
$$\langle v\,\, ,\,\, u_{i}\,\, (1 \leq i \leq r-2)\,\, ,\,\, w \rangle$$ 
forms an integral basis of $L$.

Let $h \in N^{0}$. Using $\overline{h}(\overline{u}_{i}) = \overline{u}_{i}$, 
we calculate 
$$h(v) = v\,\, ,\,\, h(u_{i}) = u_{i} + \alpha_{i}(h)v\,\, ,$$ 
where $\alpha_{i}(h)$ ($1 \leq i \leq r-2$) are integers uniquely determined 
by $h$ . 
Since ${\rm det}\, h = 1$, it follows that 
$h(w)$ is of the form: 
$$h(w) = w + \beta(h)v + \sum_{i=1}^{r-2}\gamma_{i}(h)u_{i}\,\, ,$$ 
where $\beta(h)$ and $\gamma_{i}(h)$ are also integers uniquely determined by 
$h$.  

This already shows that $N^{0}$ is unipotent. Then $N$ is of null-entropy, 
because $[N:N^{0}] < \infty$, as well. 

Let us show that $N^{0}$ is an abelian group of rank at most $r -2$. 
Varying $h$ in $N^{0}$, we can define the map $\varphi$ by: 

$$\varphi : N^{0} \longrightarrow \mathbf Z^{r-2}\,\, ;\,\, h\, \mapsto\, 
(\alpha_{i}(h))_{i=1}^{r-2}\,\, .$$

Now the next claim completes the proof of Proposition (2.9). \end{proof}

\begin{claim} \label{claim:abel} 
$\varphi$ is an injective group homomorphism. 
\end{claim} 

\begin{proof} Let $h, h' \in N^{0}$. Then by the formula above, we calculate:

$$h'h(u_{i}) = h'(u_{i} + \alpha_{i}(h)v) = h'(u_{i}) + \alpha_{i}(h)v 
= u_{i} +
 (\alpha_{i}(h') +\alpha_{i}(h))v\,\, .$$
Thus $\alpha_{i}(h'h) = \alpha_{i}(h') + \alpha_{i}(h)$ and $\varphi$ is a 
group homomorphism. 

Let us show that $\varphi$ is injective. Let $h \in {\rm Ker}\, \varphi$. Then, 
$$h(v) = v\,\, ,\,\, h(u_{i}) = u_{i}\,\, ,
\,\, h(w) = w + \beta(h)v + \sum_{i=1}^{r-2}\gamma_{i}(h)u_{i}\,\, .$$ 

It suffices to show that $h(w) = w$. 
Using $(v, u_{i}) = 0$ and $(h(x), h(y)) = (x, y)$, we calculate 
$$(w, u_{i}) = (h(w), h(u_{i})) = (w, u_{i}) + \sum_{i=1}^{r-2} 
\gamma_{j}(h)(u_{j}, u_{i})\,\, ,$$
that is, 
$$A(\gamma_{j}(h))_{j=1}^{r-2} = (0)_{j=1}^{r-2}\,\, .$$ 
Here $A := ( (u_{i}, u_{j})) \in {\rm M}(r, \mathbf Z)$. 
Since $(u_{i}, u_{j}) = (\overline{u}_{i}, \overline{u}_{j})$ 
and $\overline{L}$ is elliptic, we have $\text{det}\, A \not= 0$. Thus, 
$(\gamma_{j}(h))_{j=1}^{r-2} = (0)_{j=1}^{r-2}$, 
and therefore 
$$h(w) = w + \beta(h)v\,\, .$$ 
 
By using $(v^{2}) = 0$, we calculate 
$$(w^{2}) = (h(w)^{2}) = (w^{2}) + 2\beta(h)(w, v)\,\, .$$ 
Thus $\beta(h)(w, v) = 0$. Since $(v, v)  = (u_{i}, v) = 0$ and $L$ is non-degenerate, we have $(w, v) \not= 0$. Hence 
$\beta(h) = 0$, and therefore $h(w) = w$. 
\end{proof}

Now Theorem (2.1) follows from Proposition (2.9) applied 
for $N = G$ and $v \in \partial \mathcal P(L) \cap L \setminus \{0\}$ 
in Lemma (2.8).

\section{Groups of autmorphisms of null-entropy}

In this section, we prove Theorem (1.3) and Theorem (1.4)(1).

Let us show Theorem (1.3). Let $M$ be a projective hyperk\"ahler manifold and 
$G < {\rm Aut}\,(M)$ (resp. $G < {\rm Bir}\, (M)$) be a subgroup of null-entropy 
(resp. of null-entropy at $H^{2}$-level). Set $H := {\rm Im}(r_{NS} : G \longrightarrow {\rm O}(NS(M)))$. 
Since $M$ is projective, $\vert {\rm Ker}\,r_{NS} \vert < \infty$ by [Og2, Corollary 2.7]. 
Thus $G$ is almost abelian of rank, say $s$, iff so is $H$. However, $H$ is almost abelian of rank $\leq \rho(M) -2$ 
by Theorem (2.1).

Next we shall show the optimality of the estimate. There is a Jacobian K3 surface 
$\varphi : M \longrightarrow \mathbf P^{1}$ s.t. $\rho(M) = 20$ and the Mordell-Weil rank of $\varphi$ 
is $18$ (See e.g. [Co], [Ny], [Og1]). As we shall show in the next Section, the action of a Mordell-Weil group 
is of null-entropy. This completes the proof of Theorem (1.3). 

\begin{remark} \label{remark:np} Let $M$ be a non-projective hyperk\"ahler manifold and $G < {\rm Bir}\,(M)$. 
Then, by [Og2] (the proof of Theorem (1.2) there), we know:
\begin{list}{}{
\setlength{\leftmargin}{10pt}
\setlength{\labelwidth}{6pt}
}
\item[(1)] {\it If $NS(M)$ is elliptic, then $G$ is of null-entropy at $H^{2}$-level iff $\vert G \vert < \infty$.} 
Indeed, if $g \in G$ is of null-entropy at $H^{2}$-level, then $g^{*}\vert T(M) = id$ by [Og2, Theorem 2.4]. 
Moreover, $\vert {\rm Im}(r_{NS} : G \longrightarrow {\rm O}(NS(M)) \vert < \infty$. Thus, 
${\rm Im}(r : G \longrightarrow {\rm O}(H^{2}(M, \mathbf Z)))$ 
is finite. Then $\vert G \vert < \infty$ by [Og2, Theorem 2.3]. The other direction is clear.  
\item[(2)] {\it If $NS(M)$ is parabolic, then $G$ is always of null-entropy at $H^{2}$-level.} 
Indeed, by [Og2, Corollary 2.7], it suffices to show that $N := {\rm Im}\, (r_{NS} : G \longrightarrow {\rm O}(NS(M)))$ 
is of null-entropy (on $NS(M)$). However, in the proof of [Og2, Proposition (5.1)], we find a finite index normal unipotent 
subgroup $N^{(0)}$ of $N$. Thus, $N$ is of null-entropy. 
\end{list}
\end{remark}

Let us show Theorem (1.4)(1). Let $a(M)$ be the algebraic dimension of a K3 surface $M$. 
We argue by dividing into three cases where $NS(M)$ is elliptic, parabolic, and hyperbolic. 

If $NS(M)$ is elliptic, then the result follows from Remark (3.1)(1) above. 

If $NS(M)$ is parabolic, then a subgroup $G$ of ${\rm Aut}\,(M)$ is always of null-entropy by Remark (3.1)(2) 
above. 
On the other hand, when $NS(M)$ is parabolic, the algebraic dimension of $M$ is $1$ (by the classification 
theory of 
surfaces). The algebraic reduction map gives rise to an elliptic fibration $a : M \longrightarrow \mathbf P^{1}$ 
of $M$. This fibration is stable under ${\rm Aut}\, (M)$, and therefore so is under $G$. This completes 
the proof when 
$NS(M)$ is parabolic.

Let us consider the case where $NS(M)$ is hyperbolic. In this case $M$ is projective.

First we show "only if" part. As before, set 
$H := {\rm Im}\,(r_{NS} : G \longrightarrow {\rm O}(NS(M)))$. We may assume that $\vert G \vert = \infty$. 
Then so is $H$ by [Og2, Corollary 2.7].  Applying Theorem (2.1) for $H$, we find a primitive element 
$v \in \partial \mathcal P(M) \cap NS(M) \setminus \{0\}$ 
s.t. $h^{*}(v) = v$ for all $h \in H$. 
Here the positive cone 
$\mathcal P(M) := \mathcal P(NS(M))$ is taken to be 
the component which contains an ample class of $M$. 

By the Riemann-Roch theorem,  $v$ is represented by a non-zero 
effective divisor, say $D$, with $h^{0}(\mathcal O_{M}(D)) \geq 2$. 
Decompose $\vert D \vert$ into 
the movable part and fixed part:
$$\vert D \vert = \vert E \vert + B\,\, .$$
Then $(E^{2}) \geq 0$ and the class $[E]$ is $H$-stable. 
Thus $(E^{2}) \leq 0$ by Lemma (2.8)(1). 
Therfore $(E^{2}) = 0$. Hence $\vert E \vert$ is free and 
defines an elliptic fibration 
$$\varphi : M \longrightarrow \mathbf P^{1}$$ 
on $M$. This is $G$-stable. Moreover, $\varphi$ is the unique elliptic 
fibration satble under 
$G$. Otherwise, there is another class $[C] \in 
\partial \mathcal P(M) \cap NS(M) \setminus \{0\}$ 
such that $[C] \not\in \mathbf R_{>0}[E]$ and $H$-stable, a contradiction to 
Lemma (2.8)(2).

Next, we shall show "if part". The result is clear if $G$ is finite. 
So, we may assume 
that $\vert G \vert = \infty$ and the existence of a $G$-stable elliptic 
fibration $\varphi : M \longrightarrow \mathbf P^{1}$. 

Let $C$ be a fiber of $\varphi$. Then $[C] \in NS(M) \cap \partial \mathcal P(M) 
\setminus \{0\}$. $[C]$ is $G$-invariant as well. 
Thus, the result follows from Proposition (2.9) 
and [Og2, Corollary 2.7]. 

Let us finally show the last two statements in Theorem (1.4)(1). The second one follows from 
the contra-position of Theorem (1.3). Let us show the first one. The result is clear if $\vert G \vert < \infty$. 
So, we may assume that $\vert G \vert = \infty$. Then $G$ makes an elliptic fibration $\varphi : M \longrightarrow 
\mathbf P^{1}$ stable. Recall that any elliptic fibration on a K3 surface admits at least three singular fibers. 
(see eg. [Ca2]; This follows from $\chi_{\rm top}(M) = 24$ and $\rho(M) \geq 20$. See also [VZ] for a more 
general account.) Thus, 
${\rm Im}\,(G \longrightarrow {\rm Aut}\, (\mathbf P^{1}))$ is finite. Therefore, each orbit lies in finitely 
many fibers. This completes the proof of Theorem (1.4)(1).

\section{Autmorphism of a projective K3 surface of positive entropy}

In this section, we shall show Theorem (1.4)(2). 

Let $M$ be a K3 surface and $g \in {\rm Aut}\, (M)$.  
If $e(g) = 0$, then $\langle g \rangle$ is of null-entropy. Thus $\langle g \rangle$ has no 
Zariski dense orbit by Theorem (1.4)(1). This shows "if part" of Theorem (1.4)(2). Let us show 
"only if part" of Theorem 1.6(2). Assuming $e(g) > 0$,  we want to find a point $x \in M$ such that 
$\{g^{n}(x) \vert n \in \mathbf Z\}$ is Zariski dense in $M$. Note that $\text{ord}\ g = \infty$ y $e(g) > 0$.

Let $\mathcal F = \cup_{n \in \mathbf Z \setminus \{0\}} \mathcal F_{n}$, 
where  
$$\mathcal F_{n} := \{\,y \in M\, \vert\, g^{n}(y) = y\,\}\,\,  ,$$ 
and 
$$\mathcal C := \cup_{C \subset M, C \simeq \mathbf P^{1}}\, C\,\, .$$
Since $\text{ord}\, g = \infty$, each $\mathcal F_{n}$ is a proper closed analytic subset of $M$. 
Since $M$ is not 
uniruled, $\mathcal C$ 
is also at most countable union of $\mathbf P^{1}$. Thus, 
$\mathcal F \cup \mathcal C$ 
is a countable union of proper closed analytic subsets of $M$. Hence  
$$M \not= \mathcal F \cup \mathcal C\,\, .$$ 
Choose $x \in M \setminus (\mathcal F \cup \mathcal C)$ 
and set 
$\mathcal O(x) := \{g^{n}(x) \vert n \in \mathbf Z\}$.

The next claim will complete the proof. 

\begin{claim} \label{claim:zar} 
$\mathcal O(x)$ is Zariski dense in $M$. 
\end{claim} 

\begin{proof} Let $S$ be the Zariski closure of $\mathcal O(x)$ in $M$. 
Suppose to the contrary that $S$ is a proper subset of $M$. 
Since $\mathcal O(x)$ is an infinite set by $x \not\in \mathcal F$, 
the set $S$ is decomposed into non-empty finitely many complete irreducible curves 
and (possibly empty) finite set of closed points. Let $C$ be a $1$-dimensional 
irreducible component of $S$. Then, there is $N$ such that $g^{Nk}(C) = C$ for all 
$k \in \mathbf Z$. Note that $g^{n}(x) \in C$ for some $n$. 
(Indeed, otherwise, $S \setminus (C \cap (S\setminus C))$ 
would be a smaller closed subset containing $\mathcal O(x)$.) By the 
choice of $x$, we have $x \in g^{-n}(C)$ 
and therefore $C \simeq g^{-n}(C) \not\simeq \mathbf P^{1}$. Hence $(C^{2}) \geq 0$. 
If $(C^{2}) > 0$, then $M$ is projective and $g^{N}$, whence $g$, is of finite order 
on $NS(M)$ by Proposition (2.8)(1). Then ${\rm ord}\, g < \infty$ by 
[Og2, Corollary 2.7], a contradiction. 
If $(C^{2}) = 0$, then $C$ defines an elliptic fibration which is stable under 
$g^{N}$. Then, $g^{N}$, whence $g$, is of null-entropy by Theorem (1.4)(1), 
a contradiction. 
Therefore $S = M$. 
\end{proof}

\section{Non-abelian subgroup of automorphisms of a K3 surface}

In this section, we shall prove Theorem (1.6).

Theorem (1.6)(2) follows from Theorem 1.6(1) 
by considering the Hilbertscheme ${\rm Hilb}^{m}(S)$ of $0$-dimensional closed subschemes of length $m$ 
of a K3 surface $S$. 
Indeed, by [Be], ${\rm Hilb}^{m}(S)$ is a hyperk\"ahler manifold of dimension $2m$ having 
natural inclusions $H^{2}(S, \mathbf Z) \subset H^{2}({\rm Hilb}^{m}(S)\, , \mathbf Z)$ 
and ${\rm Aut}\, (S) \subset {\rm Aut}\, ({\rm Hilb}^{m}(S))$. Thus, ${\rm Hilb}^{m}(S)$ for a K3 surface 
$S$ in Theorem (1.6)(1) satisfies the requirement.

We shall show Theorem (1.6)(1). Let $M$ be a K3 surface and $\varphi_{i} : M \longrightarrow \mathbf P^{1}$ 
($i = 1$, $2$) be two different Jacobian fibrartions of positive Mordell-Weil rank on $M$. Note that $M$ is 
projective 
under this assumption.

Let ${\rm MW}(\varphi_{i})$ be the Mordell-Weil group of $\varphi_{i}$. 
Choose $f_{i} \in {\rm MW}(\varphi_{i})$ s.t. ${\rm ord}\, f_{i} = \infty$. 
We may regard $f_{i}$ as an element of ${\rm Aut}\, (M)$. 
Let $G := \langle f_{1}, f_{2} \rangle < {\rm Aut}\,(M)$.

First we shall show that ${\rm Aut}\, (M)$ is not almost abelian. Note that any subgroup of an almost 
abelian group 
is again almost abelian (by definition). So, it suffices to show that $G$ is not almost abelian. Suppose to 
the contrary 
that $G$ is almost abelian. Then, by definition, there are a normal subgroup $H$ of $G$ of finite index, 
a finite subgroup $N$ of $H$ and an abelian group $A$ which fit in with the exact sequence:
$$1 \longrightarrow N \longrightarrow H \longrightarrow A \longrightarrow 0\,\, .$$
Since $\vert G/H \vert < \infty$, there is a positive integer $m$ s.t. $f_{1}^{m}, f_{2}^{m} \in H$. 
Set $g_{i} := f_{i}^{m}$. Let $n := \vert N \vert$. Since $A$ is abelian, one has 
$g_{1}^{-j}g_{2}g_{1}^{j}g_{2}^{-1} \in N$ 
for each $j \in \mathbf Z$. Then, we have $n+1$-elements in $N$:
$$g_{1}^{-1}g_{2}g_{1}g_{2}^{-1}\,\, ,\,\, g_{1}^{-2}g_{2}g_{1}^{2}g_{2}^{-1}\,\, ,\,\, \cdots \,\, ,\,\, 
g_{1}^{-n}g_{2}g_{1}^{n}g_{2}^{-1}\,\, ,\,\, g_{1}^{-n-1}g_{2}g_{1}^{n+1}g_{2}^{-1}\,\, \in N\,\, .$$ 
Thus, at least two of them have to coinside. Hence $g_{1}^{k}g_{2} = g_{2}g_{1}^{k}$ for some positive integer $k$. 
Set $g := g_{1}^{k}$ and $h := g_{2}$. Then $gh = hg$ and $g^{*}(e_{1}) = e_{1}$, $h^{*}(e_{2}) = e_{2}$. 
Here $e_{i}$ is the class of a general fiber of $\varphi_{i}$ ($i = 1$, $2$). Set $e_{3} := h^{*}(e_{1})$. 
Then $e_{3}$ is a (primitive) class of elliptic pencils on $M$. If $e_{3} = e_{1}$, 
then $e_{1}$ and $e_{2}$ are both $h$-stable. Therefore, $e_{1} + e_{2}$ 
is also $h$-stable. However, since $(e_{1} + e_{2})^{2} > 0$, we have then ${\rm ord}\, h < \infty$ 
on $NS(M)$, whence ${\rm ord}\, g_{2} < \infty$ by [Og2, Corollary 2.7], a contradiction to 
${\rm ord}\, g_{2} = \infty$. 
Consider next the case $e_{3} \not= e_{1}$, i.e. the case where $e_{1}$ and $e_{3}$ correspond different elliptic 
pencils. Using $g^{*}h^{*} = h^{*}g^{*}$ 
and $g^{*}(e_{1}) = e_{1}$, 
one has $g^{*}(e_{3}) = e_{3}$. Then $e_{1} + e_{3}$ is $g$-stable and $(e_{1} + e_{3})^{2} > 0$, 
a contradiction for the same reason as above. Thus, $G$ is not almost abelian.

Since $G$ is not almost abelian, $G$ is of positive entropy by the contra-position of Theorem (1.3). 

Next we shall show that a K3 surface $M$ of Picard number $20$ has at least two different Jacobian fibrations 
of positive Mordell-Weil rank. Note that $M$ is necessarily projective. 
In what follows, $U$ stands for an even unimodular hyperbolic lattice of rank 
$2$, and $A_{*}$, $D_{*}$, $E_{*}$ stands for the negative definite root 
lattice whose basis is given by the corresponding Dynkin diagram. Since the transcendental lattice $T(M)$ is 
positive definite and of rank $2$, 
$T(M)$ can be primitively embedded into $U^{\oplus 2}$, and $T(M)$ can be also uniquely primitively embedded into 
the K3 lattice $\Lambda_{\rm K3} := U^{3} \oplus E_{8}(-1)^{2}$ (See e.g. [Ni]). Note that $H^{2}(M, \mathbf Z) \simeq \Lambda_{\rm K3}$.  
Thus the N\'eron-Severi lattice $NS(M)$ of $M$ is of the form: 
$$NS(M) = U \oplus E_{8}(-1)^{\oplus 2} \oplus N\,\, ,$$ 
where $N$ is a negative definite lattice of rank $2$. Put $d := {\rm det}\, N$. If $d = 3, 4$, then 
$N = A_{2}$, $A_{1}^{2}$, and if $d \not= 3$, $4$, then $N$ is not a root lattice. Recall that the 
$\mathbf Q$-rational hull 
of the ample cone of a projective K3 surface is the fundamental domain of the group of $(-2)$-reflections 
on the rational hull of the positive cone (See e.g. [BPV]). Using this fact and the Riemann-Roch formula, one finds 
a Jacobian fibration $\varphi_{1} : M \longrightarrow \mathbf P^{1}$ whose reducible singular fibers are $II^{*} + II^{*}$ if $d \not= 
3, 4$, $II^{*} + II^{*} + I_{3}$ or $II^{*} + II^{*} + IV$ if $d = 3$ and $II^{*} + II^{*} + A + B$ 
($A, B \in \{I_{2}, III\}$) if $d = 4$. Let ${\rm MW}\,(\varphi_{1})$ be the Mordell-Weil group of $\varphi_{1}$. Then 
${\rm rank}\, {\rm MW}\,(\varphi_{1}) > 0$ unless $d = 3$, $4$ 
by Shioda's rank formula (See for instance [Sh]). Note also that this $\varphi_{1}$ is essentially the same Jacobian 
fibration studied in [SI].

Let us find another Jacobian fibrations of positive Mordell-Weil rank. Join two $II^{*}$ singular fibers 
of $\varphi_{1}$ by the $0$-section and throw out the components of multiplicity 
$2$ at the edge of two $II^{*}$, say $C_{1}$ in one $II^{*}$ and $C_{2}$ in the other $II^{*}$. 
In this way, one obtains a divisor of Kodaira's type $I_{12}^{*}$, say, $D$, on $M$.  
The pencil $\vert D \vert$ gives rise to a Jacobian fibration 
$\varphi_{2} : M \longrightarrow \mathbf P^{1}$ with reducible singular fiber of type $I_{12}^{*}$ and with sections 
$C_{1}$, $C_{2}$. Take $C_{1}$ as the $0$ and consider $C_{2}$ as an element $P$ of ${\rm MW}\,(\varphi_{2})$. We have 
that ${\rm rank}\, {\rm MW}\, (\varphi_{2}) > 0$ (regardless of $d$).

Indeed, otherwise, by using Shioda's rank formula, we see that the remaining reducible singular fibers 
of $\varphi_{2}$ are either 
$I_{3}$, $VI$, $I_{2} + I_{2}$, $I_{2} + III$, or $III + III$. In each case, the value  
$\langle P , P \rangle$ of Shioda's height pairing of $P$ is positive, as it is checked by his explicit height pairing 
formula. This, however, implies 
${\rm ord}\, P = \infty$, a contradiction to the assumption. (See [Sh] for the definition and basic properties on 
Shioda's height pairing on the Mordell-Weil group.)

Thus we obtain two Jacobian fibrations $\varphi_{1}$, $\varphi_{2}$ of positive Mordell-Weil rank 
unless $d \not= 3$, $4$. When $d = 3$, $4$, Vinberg [Vi] calculated the full automorphism group ${\rm Aut}\, (M)$. It 
contains a subgroup $H$ isomorphic to 
the free product $\mathbf Z_{2} * \mathbf Z_{2} * \mathbf Z_{2}$ of three cyclic groups of order $2$. 
This $H$ is far from being almost abelian. Thus, by Theorem (1.4)(1), $H$ does 
not make the Jacobian fibration $\varphi_{2}$ stable. Now, by transforming $\varphi_{2}$ by 
$H$, one can find another Jacobian fibration of positive rank also when $d = 3$, $4$. This completes the proof 
of Theorem (1.6)(1). 

\begin{remark} \label{remark:cn} Under the same asumption in Theorem (1.6)(1), 
a K3 surface $M$ has a dense ${\rm Aut}\,(M)$ orbit in the Euclidean topology. 
This remarkable result has been found by 
Cantat [Ca2] in a slightly more general situation. 
\end{remark}


\end{document}